\documentclass[10pt]{article}
\newtheorem{thm}{Theorem} %[section]

\newcommand{\QED}{\hfill $\fbox{}$\vspace*{3mm}}
\newcommand{\Za}{\hbox{\bf Z}_a}
\newcommand{\N}{\hbox{\bf N}}
\newcommand{\Z}{\hbox{\bf Z}}
\newcommand{\G}{\mbox{$\Gamma$}}
\newcommand{\subgp}[1]{\langle{#1}\rangle}
\newcommand{\beeq}{\begin{eqnarray*}}
\newcommand{\eneq}{\end{eqnarray*}}
\newcommand{\proof}{\noindent {\it Proof.\hspace{4mm}}}
\newcommand{\example}{\noindent{\bf Example.\hspace{4mm}}}
\newcommand{\qfd}{\hfill $\fbox{}$\vspace{4mm}}\def\newpic#1{%
\def\emline##1##2##3##4##5##6{%
\put(##1,##2){\special{em:point #1##3}}%
\put(##4,##5){\special{em:point #1##6}}%
\special{em:line #1##3,#1##6}}}
\newpic{}
\def\emline#1#2#3#4#5#6{%
\put(#1,#2){\special{em:moveto}}%
\put(#4,#5){\special{em:lineto}}}
\def\newpic#1{}
\title{On a $\vec{C}_4$-ultrahomogeneous oriented graph}
\author{Italo J. Dejter
\\ University of Puerto Rico \\ Rio Piedras, PR 00931-3355 \\ idejter@uprrp.edu}
\date{}

\begin{document}
\maketitle

\begin{abstract}
The notion of a $\mathcal C$-ul\-tra\-ho\-mo\-ge\-neous graph, due to Isaksen et al., is adapted for digraphs, and then a strongly connected $\vec{C}_4$-ul\-tra\-ho\-mo\-ge\-neous oriented graph on 168 vertices and 126 pairwise arc-disjoint 4-cycles is presented, with regular indegree and outdegree 3 and no circuits of lengths 2 and 3, by altering a definition of the Coxeter graph via pencils of ordered lines of the Fano plane in which pencils are replaced by ordered pencils.
\end{abstract}


\section{Introduction}

The study of ultrahomogeneous graphs (resp. digraphs) can be traced back to \cite{Sheh}, \cite{Gard}, \cite{Ronse} and, \cite{GK}, (resp. \cite{Fra}, \cite{Lach} and \cite{Cher}). In \cite{I}, $\mathcal C$-ultrahomogeneous graphs are defined and subsequently treated when ${\mathcal C}=$ collection of either {\bf(a)} complete graphs, or {\bf(b)} disjoint unions of complete graphs, or {\bf(c)} complements of those unions.
In \cite{D1}, a $\{K_4,K_{2,2,2}\}$-ultrahomogeneous graph on 42 vertices, 42 copies of $K_4$ and 21 copies of $K_{2,2,2}$ is given that fastens objects of (a) and (c), namely $K_4$ and $K_{2,2,2}$, respectively, over copies of $K_2$.

In the present note and in \cite{DD}, the notion of a $\mathcal C$-ultrahomogeneous graph is extended as follows:
Given a collection $\mathcal C$ of (di)graphs closed under isomorphisms, a (di)graph $G$ is $\mathcal C$-{\it ul\-tra\-ho\-mo\-ge\-neous} (or $\mathcal C$-UH) if every isomorphism between two $G$-induced members of $\mathcal C$ extends to an auto\-mor\-phism of $G$. If ${\mathcal C}=\{H\}$ is the isomorphism class of a (di)graph $H$, such a $G$ is said to be $\{H\}$-UH or $H$-UH.

In \cite{DD}, the cubic distance-transitive graphs are shown to be $C_g$-UH graphs, where $C_g$ stands for cycle of minimum length, i.e. realizing the girth $g$; moreover, all these graphs but for the Petersen, Heawood and Foster graphs are shown to be $\vec{C}_g$-UH digraphs,
which allows the construction of novel $\mathcal C$-UH graphs, in continuation to the work of \cite{D1}, including a $\{K_4,L(Q_3)\}$-UH graph on 102 vertices that fastens 102 copies of $K_4$ and 102 copies of the cuboctahedral graph $L(Q_3)$ over copies of $K_3$, obtained from the Biggs-Smith graph, (\cite{Bi}), by unzipping, powering and zipping back a collection of oriented g-cycles provided by the initial results. However, these graphs are undirected, so they are not properly digraphs.

In this note, a presentation of the Coxeter graph $Cox$ is modified
to provide a strongly connected $\vec{C}_4$-UH oriented graph $D$ on 168 vertices, 126 pairwise arc-disjoint 4-cycles, with regular indegree and outdegree 3. In contrast, the construction of \cite{D1} used ordered pencils of unordered lines, instead.

We take the Fano plane $\mathcal F$ as having point set $J_7=\Z_7$ (the cyclic group mod 7) and point-line correspondence $\phi(j)=\{(j+1),(j+2),(j+4)\}$, for every $j\in\Z_7$
in order to color the vertices and edges of $Cox$ as in Figure 1.
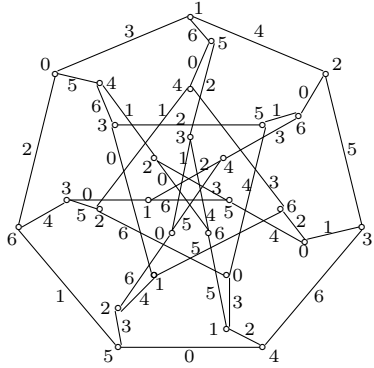
\begin{figure}[htp]
\unitlength=0.40mm
\special{em:linewidth 0.4pt}
\linethickness{0.4pt}
\begin{picture}(213.00,117.00)
\put(142.00,28.00){\circle{2.00}}
\put(166.00,28.00){\circle{2.00}}
\put(124.00,50.00){\circle{2.00}}
\put(178.00,78.00){\circle{2.00}}
\put(142.00,67.00){\circle{2.00}}
\put(154.00,90.00){\circle{2.00}}
\put(160.00,42.00){\circle{2.00}}
\put(148.00,42.00){\circle{2.00}}
\put(113.00,53.00){\circle{2.00}}
\put(184.00,50.00){\circle{2.00}}
\put(165.00,67.00){\circle{2.00}}
\put(154.00,74.00){\circle{2.00}}
\put(142.00,28.00){\circle{2.00}}
\put(109.00,95.00){\circle{2.00}}
\put(154.00,114.00){\circle{2.00}}
\put(199.00,95.00){\circle{2.00}}
\put(97.00,44.00){\circle{2.00}}
\put(211.00,44.00){\circle{2.00}}
\put(130.00,4.00){\circle{2.00}}
\put(178.00,4.00){\circle{2.00}}
\emline{155.00}{114.00}{1}{198.00}{96.00}{2}
\emline{110.00}{96.00}{3}{153.00}{114.00}{4}
\emline{199.00}{94.00}{5}{211.00}{45.00}{6}
\emline{109.00}{94.00}{7}{97.00}{45.00}{8}
\emline{97.00}{43.00}{9}{129.00}{5.00}{10}
\emline{131.00}{4.00}{11}{177.00}{4.00}{12}
\emline{179.00}{5.00}{13}{211.00}{43.00}{14}
\put(129.00,78.00){\circle{2.00}}
\emline{155.00}{90.00}{15}{185.00}{51.00}{16}
\emline{184.00}{49.00}{17}{143.00}{28.00}{18}
\emline{141.00}{29.00}{19}{128.00}{77.00}{20}
\emline{130.00}{78.00}{21}{177.00}{78.00}{22}
\emline{179.00}{77.00}{23}{167.00}{28.00}{24}
\emline{165.00}{28.00}{25}{124.00}{49.00}{26}
\emline{123.00}{51.00}{27}{153.00}{90.00}{28}
\put(167.00,53.00){\circle{2.00}}
\emline{154.00}{73.00}{29}{160.00}{43.00}{30}
\emline{160.00}{43.00}{31}{143.00}{66.00}{32}
\emline{168.00}{52.00}{33}{191.00}{39.00}{34}
\emline{139.00}{53.00}{35}{114.00}{53.00}{36}
\emline{141.00}{53.00}{37}{164.00}{66.00}{38}
\emline{164.00}{66.00}{39}{148.00}{43.00}{40}
\emline{148.00}{43.00}{41}{154.00}{73.00}{42}
\put(161.00,106.00){\circle{2.00}}
\put(190.00,81.00){\circle{2.00}}
\emline{166.00}{67.00}{43}{189.00}{80.00}{44}
\put(192.00,39.00){\circle{2.00}}
\emline{143.00}{66.00}{45}{166.00}{53.00}{46}
\emline{160.00}{41.00}{47}{166.00}{11.00}{48}
\put(166.00,10.00){\circle{2.00}}
\put(130.00,17.00){\circle{2.00}}
\emline{147.00}{41.00}{49}{131.00}{18.00}{50}
\put(140.00,53.00){\circle{2.00}}
\emline{166.00}{53.00}{51}{141.00}{53.00}{52}
\put(124.00,92.00){\circle{2.00}}
\emline{142.00}{68.00}{53}{125.00}{91.00}{54}
\emline{154.00}{113.00}{55}{160.00}{106.00}{56}
\emline{198.00}{94.00}{57}{191.00}{82.00}{58}
\emline{210.00}{44.00}{59}{193.00}{40.00}{60}
\emline{177.00}{5.00}{61}{167.00}{10.00}{62}
\emline{131.00}{5.00}{63}{129.00}{16.00}{64}
\emline{98.00}{44.00}{65}{112.00}{52.00}{66}
\emline{110.00}{95.00}{67}{123.00}{92.00}{68}
\emline{154.00}{75.00}{69}{162.00}{105.00}{70}
\emline{154.00}{91.00}{71}{160.00}{105.00}{72}
\emline{179.00}{79.00}{73}{189.00}{82.00}{74}
\emline{185.00}{49.00}{75}{192.00}{40.00}{76}
\emline{167.00}{27.00}{77}{167.00}{11.00}{78}
\emline{131.00}{16.00}{79}{142.00}{27.00}{80}
\emline{123.00}{91.00}{81}{128.00}{79.00}{82}
\put(157.00,117.00){\makebox(0,0)[cc]{$_1$}}
\put(177.00,109.00){\makebox(0,0)[cc]{$_4$}}
\put(203.00,98.00){\makebox(0,0)[cc]{$_2$}}
\put(213.00,40.00){\makebox(0,0)[cc]{$_3$}}
\put(182.00,1.00){\makebox(0,0)[cc]{$_4$}}
\put(197.00,21.00){\makebox(0,0)[cc]{$_6$}}
\put(208.00,70.00){\makebox(0,0)[cc]{$_5$}}
\put(154.00,1.00){\makebox(0,0)[cc]{$_0$}}
\put(111.00,21.00){\makebox(0,0)[cc]{$_1$}}
\put(95.00,40.00){\makebox(0,0)[cc]{$_6$}}
\put(100.00,70.00){\makebox(0,0)[cc]{$_2$}}
\put(106.00,98.00){\makebox(0,0)[cc]{$_0$}}
\put(134.00,109.00){\makebox(0,0)[cc]{$_3$}}
\put(155.00,108.00){\makebox(0,0)[cc]{$_6$}}
\put(165.00,105.00){\makebox(0,0)[cc]{$_5$}}
\put(150.00,91.00){\makebox(0,0)[cc]{$_4$}}
\put(155.00,99.00){\makebox(0,0)[cc]{$_0$}}
\put(161.00,91.00){\makebox(0,0)[cc]{$_2$}}
\put(151.00,74.00){\makebox(0,0)[cc]{$_3$}}
\put(167.00,64.00){\makebox(0,0)[cc]{$_4$}}
\put(167.00,49.00){\makebox(0,0)[cc]{$_5$}}
\put(164.00,42.00){\makebox(0,0)[cc]{$_6$}}
\put(144.00,42.00){\makebox(0,0)[cc]{$_0$}}
\put(140.00,49.00){\makebox(0,0)[cc]{$_1$}}
\put(140.00,64.00){\makebox(0,0)[cc]{$_2$}}
\put(152.00,67.00){\makebox(0,0)[cc]{$_1$}}
\put(145.00,60.00){\makebox(0,0)[cc]{$_0$}}
\put(146.00,51.00){\makebox(0,0)[cc]{$_6$}}
\put(153.00,45.00){\makebox(0,0)[cc]{$_5$}}
\put(161.00,48.00){\makebox(0,0)[cc]{$_4$}}
\put(163.00,57.00){\makebox(0,0)[cc]{$_3$}}
\put(159.00,65.00){\makebox(0,0)[cc]{$_2$}}
\put(183.00,83.00){\makebox(0,0)[cc]{$_1$}}
\put(177.00,82.00){\makebox(0,0)[cc]{$_5$}}
\put(191.00,77.00){\makebox(0,0)[cc]{$_6$}}
\put(192.00,89.00){\makebox(0,0)[cc]{$_0$}}
\put(134.00,83.00){\makebox(0,0)[cc]{$_1$}}
\put(151.00,80.00){\makebox(0,0)[cc]{$_2$}}
\put(128.00,92.00){\makebox(0,0)[cc]{$_4$}}
\put(123.00,84.00){\makebox(0,0)[cc]{$_6$}}
\put(125.00,78.00){\makebox(0,0)[cc]{$_3$}}
\put(128.00,67.00){\makebox(0,0)[cc]{$_0$}}
\put(144.00,25.00){\makebox(0,0)[cc]{$_1$}}
\put(126.00,17.00){\makebox(0,0)[cc]{$_2$}}
\put(134.00,27.00){\makebox(0,0)[cc]{$_6$}}
\put(133.00,11.00){\makebox(0,0)[cc]{$_3$}}
\put(162.00,10.00){\makebox(0,0)[cc]{$_1$}}
\put(174.00,10.00){\makebox(0,0)[cc]{$_2$}}
\put(192.00,36.00){\makebox(0,0)[cc]{$_0$}}
\put(200.00,44.00){\makebox(0,0)[cc]{$_1$}}
\put(191.00,46.00){\makebox(0,0)[cc]{$_2$}}
\put(188.00,51.00){\makebox(0,0)[cc]{$_6$}}
\put(182.00,59.00){\makebox(0,0)[cc]{$_3$}}
\put(124.00,46.00){\makebox(0,0)[cc]{$_2$}}
\put(132.00,42.00){\makebox(0,0)[cc]{$_6$}}
\put(118.00,48.00){\makebox(0,0)[cc]{$_5$}}
\put(113.00,57.00){\makebox(0,0)[cc]{$_3$}}
\put(120.00,55.00){\makebox(0,0)[cc]{$_0$}}
\put(182.00,41.00){\makebox(0,0)[cc]{$_4$}}
\put(145.00,83.00){\makebox(0,0)[cc]{$_1$}}
\put(139.00,20.00){\makebox(0,0)[cc]{$_4$}}
\put(170.00,19.00){\makebox(0,0)[cc]{$_3$}}
\put(170.00,28.00){\makebox(0,0)[cc]{$_0$}}
\put(156.00,37.00){\makebox(0,0)[cc]{$_5$}}
\put(115.00,91.00){\makebox(0,0)[cc]{$_5$}}
\put(107.00,46.00){\makebox(0,0)[cc]{$_4$}}
\put(184.00,74.00){\makebox(0,0)[cc]{$_3$}}
\put(161.00,22.00){\makebox(0,0)[cc]{$_5$}}
\put(173.00,58.00){\makebox(0,0)[cc]{$_4$}}
\put(127.00,1.00){\makebox(0,0)[cc]{$_5$}}
\emline{114.00}{53.00}{83}{123.00}{50.00}{84}
\end{picture}
\caption{Coloring the vertices and edges of $Cox$ with elements of $\mathcal F$}
\end{figure}

This figure shows that each vertex $v$ of $Cox$ can be considered as an unordered pencil of ordered lines of $\mathcal F$, (brackets and commas avoided now):
\begin{equation}xb_1c_1,\,\,\,xb_2c_2,\,\,\,xb_0c_0,\label{p}\end{equation} corresponding to the three edges $e_1,e_2,e_0$ incident to $v$, respectively, and denoted by $[x,b_1c_1,b_2c_2,b_0c_0]$, where $x$ is the color of $v$ in the figure, with $b_i$ and $c_i$ as the colors of the edge $e_i$ and the endvertex of $e_i$ other than $v$, for $i\in\{1,2,0\}$.

Moreover, two such vertices are adjacent in $Cox$ if they can be written $[x,b_1c_1,b_2c_2,b_0c_0]$ and $[x',b'_1c'_1,b'_2c'_2,b'_0c'_0]$ (perhaps by means of a permutation of the entries $b_ic_i$) in such a way that $\{b_i,c_i\}\cap\{b'_i,c'_i\}$ is constituted by just one element $d_i$, for each $i\in\{1,2,0\}$, and the resulting triple $d_1d_2d_0$ is a line of $\mathcal F$.

\section{Presentation of a $\vec{C_4}$-UH digraph}

Consider the oriented graph $D$ whose vertices are the {\it ordered} pencils
of ordered lines of $\mathcal F$, as in (1) above. Each such vertex will be denoted  $(x,b_1c_1,b_2c_2,b_0c_0),$ where $b_1b_2b_0$
is a line of $\mathcal F$. An arc between two vertices of $D$, say from  $(x,b_1c_1,b_2c_2,b_0c_0)$ and $(x',b'_1c'_1,b'_2c'_2,b'_0c'_0)$, is established if and only if
$$\begin{array}{rlll}
x=c'_i, & b'_{i+1}=c_{i+1}, & b'_{i-1}=c_{i-1}, & b'_i=b_i, \\
x'=c_i, & c'_{i+1}=b_{i-1}, & c'_{i-1}=b_{i+1}, &
 \end{array}$$ for some, $i\in\{1,2,0\}$. This way,
we obtain oriented 4-cycles in $D$, such as $$((0,26,54,31),(6,20,43,15),(0,26,31,54),(6,20,15,43)).$$
A simplified notation for the vertices $(x,yz,uv,pq)$ of $D$ is $yup_x$. With such a notation,
the adjacency sub-list of $D$ departing from the vertices of the form $yup_0$ is
(with rows indicated $a,b,c,d,e,f$, to be used below):
$$\begin{array}{cccc}
^{124_0\,:\,165_3, 325_6, 364_5;}_{142_0\,:\,156_3, 346_5, 352_6;} &
^{235_0\,:\,214_6, 634_1, 615_6;}_{253_0\,:\,241_6, 651_4, 643_6;} &
^{346_0\,:\,352_1, 142_5, 156_2;}_{364_0\,:\,325_1, 165_2, 124_5;} &
^{156_0\,:\,142_3, 352_4, 346_2;}_{165_0\,:\,124_3, 364_2, 325_4;}\vspace*{0.6mm} \\
^{214_0\,:\,235_6, 615_3, 634_5;}_{241_0\,:\,253_6, 643_5, 651_3;} &
^{325_0\,:\,364_1, 124_6, 165_1;}_{352_0\,:\,346_1, 156_4, 142_1;} &
^{436_0\,:\,412_5, 532_1, 516_2;}_{463_0\,:\,421_5, 561_2, 523_1;} &
^{516_0\,:\,532_4, 412_3, 436_2;}_{561_0\,:\,523_4, 463_2, 421_3;}\vspace*{0.6mm} \\
^{412_0\,:\,436_5, 516_3, 532_6;}_{421_0\,:\,463_5, 523_6, 561_3;} &
^{523_0\,:\,561_4, 421_6, 463_4;}_{532_0\,:\,516_4, 436_1, 412_4;} &
^{634_0\,:\,615_2, 235_1, 214_5;}_{643_0\,:\,651_2, 241_5, 253_1;} &
^{615_0\,:\,634_2, 214_3, 235_4;}_{651_0\,:\,643_2, 253_4, 241_3.}
\end{array}$$
From this sub-list, the adjacency list of $D$, for its $168=24\times 7$ vertices, is obtained via translations mod 7. Let us represent each vertex $yup_0$ of $D$ by means of a symbol $j_i$, where $i=a,b,c,d,e,f$ stands for the successive rows of the table above and
$j=\phi^{-1}(yup)\in\{0,1,2,4\}$. These symbols $j_i$ are assigned to the lines $yup$ avoiding $0\in{\mathcal F}$, and thus to the vertices $yup_0$, as follows:

$$\begin{array}{l|ccccc}
_{j_i}&_{j=0}&_{j=1}&_{j=2}&_{j=4} \\ \hline
&&&& \\
^{i=a}_{i=b}&^{124}_{142}&^{235}_{253}&^{346}_{364}&^{156}_{165}\\
^{i=c}_{i=d}&^{214}_{241}&^{325}_{352}&^{436}_{463}&^{516}_{561}\\
^{i=e}_{i=f}&^{412}_{421}&^{523}_{532}&^{634}_{643}&^{615}_{651}
\end{array}$$
\begin{figure}[htp]
\unitlength=0.60mm
\special{em:linewidth 0.4pt}
\linethickness{0.4pt}
\begin{picture}(198.00,57.00)
\put(5.00,4.00){\circle{2.00}}
\put(35.00,4.00){\circle{2.00}}
\put(65.00,4.00){\circle{2.00}}
\put(95.00,4.00){\circle{2.00}}
\put(5.00,24.00){\circle{2.00}}
\put(35.00,24.00){\circle{2.00}}
\put(65.00,24.00){\circle{2.00}}
\put(95.00,24.00){\circle{2.00}}
\emline{5.00}{5.00}{1}{5.00}{23.00}{2}
\emline{35.00}{5.00}{3}{35.00}{23.00}{4}
\put(25.00,14.00){\circle{2.00}}
\put(45.00,14.00){\circle{2.00}}
\emline{26.00}{15.00}{5}{34.00}{23.00}{6}
\emline{36.00}{23.00}{7}{44.00}{15.00}{8}
\emline{44.00}{13.00}{9}{36.00}{5.00}{10}
\emline{34.00}{5.00}{11}{26.00}{13.00}{12}
\put(55.00,14.00){\circle{2.00}}
\put(75.00,14.00){\circle{2.00}}
\emline{56.00}{15.00}{13}{64.00}{23.00}{14}
\emline{66.00}{23.00}{15}{74.00}{15.00}{16}
\emline{74.00}{13.00}{17}{66.00}{5.00}{18}
\emline{64.00}{5.00}{19}{56.00}{13.00}{20}
\put(85.00,14.00){\circle{2.00}}
\emline{86.00}{15.00}{21}{94.00}{23.00}{22}
\emline{94.00}{5.00}{23}{86.00}{13.00}{24}
\put(15.00,14.00){\circle{2.00}}
\emline{6.00}{23.00}{25}{14.00}{15.00}{26}
\emline{14.00}{13.00}{27}{6.00}{5.00}{28}
\put(5.00,34.00){\circle{2.00}}
\put(35.00,34.00){\circle{2.00}}
\put(65.00,34.00){\circle{2.00}}
\put(95.00,34.00){\circle{2.00}}
\put(5.00,54.00){\circle{2.00}}
\put(65.00,54.00){\circle{2.00}}
\put(95.00,54.00){\circle{2.00}}
\emline{5.00}{35.00}{29}{5.00}{53.00}{30}
\emline{35.00}{35.00}{31}{35.00}{53.00}{32}
\put(25.00,44.00){\circle{2.00}}
\put(45.00,44.00){\circle{2.00}}
\put(55.00,44.00){\circle{2.00}}
\put(75.00,44.00){\circle{2.00}}
\put(85.00,44.00){\circle{2.00}}
\put(15.00,44.00){\circle{2.00}}
\put(5.00,57.00){\makebox(0,0)[cc]{$_{1_b}$}}
\put(35.00,57.00){\makebox(0,0)[cc]{$_{1_c}$}}
\put(65.00,57.00){\makebox(0,0)[cc]{$_{1_f}$}}
\put(95.00,57.00){\makebox(0,0)[cc]{$_{1_b}$}}
\put(5.00,31.00){\makebox(0,0)[cc]{$_{1_d}$}}
\put(35.00,31.00){\makebox(0,0)[cc]{$_{1_a}$}}
\put(65.00,31.00){\makebox(0,0)[cc]{$_{1_e}$}}
\put(95.00,31.00){\makebox(0,0)[cc]{$_{1_d}$}}
\put(7.00,43.00){\makebox(0,0)[cc]{$_6$}}
\put(10.00,42.00){\makebox(0,0)[cc]{$_4$}}
\put(10.00,47.00){\makebox(0,0)[cc]{$_3$}}
\put(37.00,43.00){\makebox(0,0)[cc]{$_6$}}
\put(40.00,42.00){\makebox(0,0)[cc]{$_4$}}
\put(40.00,47.00){\makebox(0,0)[cc]{$_3$}}
\put(33.00,45.00){\makebox(0,0)[cc]{$_6$}}
\put(31.00,41.00){\makebox(0,0)[cc]{$_3$}}
\put(31.00,47.00){\makebox(0,0)[cc]{$_4$}}
\put(67.00,43.00){\makebox(0,0)[cc]{$_6$}}
\put(70.00,42.00){\makebox(0,0)[cc]{$_4$}}
\put(70.00,47.00){\makebox(0,0)[cc]{$_3$}}
\put(63.00,45.00){\makebox(0,0)[cc]{$_6$}}
\put(61.00,41.00){\makebox(0,0)[cc]{$_3$}}
\put(61.00,47.00){\makebox(0,0)[cc]{$_4$}}
\put(93.00,45.00){\makebox(0,0)[cc]{$_6$}}
\put(91.00,41.00){\makebox(0,0)[cc]{$_3$}}
\put(91.00,47.00){\makebox(0,0)[cc]{$_4$}}
\put(18.00,41.00){\makebox(0,0)[cc]{$_{0_e}$}}
\put(24.00,41.00){\makebox(0,0)[cc]{$_{0_f}$}}
\put(48.00,41.00){\makebox(0,0)[cc]{$_{0_d}$}}
\put(54.00,41.00){\makebox(0,0)[cc]{$_{0_b}$}}
\put(78.00,41.00){\makebox(0,0)[cc]{$_{0_a}$}}
\put(84.00,41.00){\makebox(0,0)[cc]{$_{0_c}$}}
\put(20.00,52.00){\makebox(0,0)[cc]{$_1$}}
\put(20.00,36.00){\makebox(0,0)[cc]{$_1$}}
\put(50.00,52.00){\makebox(0,0)[cc]{$_1$}}
\put(50.00,36.00){\makebox(0,0)[cc]{$_1$}}
\put(80.00,52.00){\makebox(0,0)[cc]{$_1$}}
\put(80.00,36.00){\makebox(0,0)[cc]{$_1$}}
\put(5.00,27.00){\makebox(0,0)[cc]{$_{2_c}$}}
\put(35.00,27.00){\makebox(0,0)[cc]{$_{2_f}$}}
\put(65.00,27.00){\makebox(0,0)[cc]{$_{2_b}$}}
\put(95.00,27.00){\makebox(0,0)[cc]{$_{2_c}$}}
\put(5.00,1.00){\makebox(0,0)[cc]{$_{2_e}$}}
\put(35.00,1.00){\makebox(0,0)[cc]{$_{2_d}$}}
\put(65.00,1.00){\makebox(0,0)[cc]{$_{2_a}$}}
\put(95.00,1.00){\makebox(0,0)[cc]{$_{2_e}$}}
\put(10.00,12.00){\makebox(0,0)[cc]{$_1$}}
\put(10.00,17.00){\makebox(0,0)[cc]{$_6$}}
\put(40.00,12.00){\makebox(0,0)[cc]{$_1$}}
\put(40.00,17.00){\makebox(0,0)[cc]{$_6$}}
\put(33.00,15.00){\makebox(0,0)[cc]{$_5$}}
\put(31.00,11.00){\makebox(0,0)[cc]{$_6$}}
\put(31.00,17.00){\makebox(0,0)[cc]{$_1$}}
\put(70.00,12.00){\makebox(0,0)[cc]{$_1$}}
\put(70.00,17.00){\makebox(0,0)[cc]{$_6$}}
\put(63.00,15.00){\makebox(0,0)[cc]{$_5$}}
\put(61.00,11.00){\makebox(0,0)[cc]{$_6$}}
\put(61.00,17.00){\makebox(0,0)[cc]{$_1$}}
\put(93.00,15.00){\makebox(0,0)[cc]{$_5$}}
\put(91.00,11.00){\makebox(0,0)[cc]{$_6$}}
\put(91.00,17.00){\makebox(0,0)[cc]{$_1$}}
\put(18.00,11.00){\makebox(0,0)[cc]{$_{0_a}$}}
\put(24.00,11.00){\makebox(0,0)[cc]{$_{0_b}$}}
\put(48.00,11.00){\makebox(0,0)[cc]{$_{0_e}$}}
\put(54.00,11.00){\makebox(0,0)[cc]{$_{0_c}$}}
\put(78.00,11.00){\makebox(0,0)[cc]{$_{0_d}$}}
\put(84.00,11.00){\makebox(0,0)[cc]{$_{0_f}$}}
\put(20.00,22.00){\makebox(0,0)[cc]{$_2$}}
\put(20.00,6.00){\makebox(0,0)[cc]{$_2$}}
\put(50.00,22.00){\makebox(0,0)[cc]{$_2$}}
\put(50.00,6.00){\makebox(0,0)[cc]{$_2$}}
\put(80.00,22.00){\makebox(0,0)[cc]{$_2$}}
\put(80.00,6.00){\makebox(0,0)[cc]{$_2$}}
\put(7.00,13.00){\makebox(0,0)[cc]{$_5$}}
\put(37.00,13.00){\makebox(0,0)[cc]{$_5$}}
\put(67.00,13.00){\makebox(0,0)[cc]{$_5$}}
\put(107.00,19.00){\circle{2.00}}
\put(137.00,19.00){\circle{2.00}}
\put(167.00,19.00){\circle{2.00}}
\put(197.00,19.00){\circle{2.00}}
\put(107.00,39.00){\circle{2.00}}
\put(137.00,39.00){\circle{2.00}}
\put(167.00,39.00){\circle{2.00}}
\put(197.00,39.00){\circle{2.00}}
\put(127.00,29.00){\circle{2.00}}
\put(147.00,29.00){\circle{2.00}}
\put(157.00,29.00){\circle{2.00}}
\put(177.00,29.00){\circle{2.00}}
\put(187.00,29.00){\circle{2.00}}
\put(117.00,29.00){\circle{2.00}}
\put(107.00,42.00){\makebox(0,0)[cc]{$_{4_b}$}}
\put(137.00,42.00){\makebox(0,0)[cc]{$_{4_c}$}}
\put(167.00,42.00){\makebox(0,0)[cc]{$_{4_f}$}}
\put(197.00,42.00){\makebox(0,0)[cc]{$_{4_b}$}}
\put(107.00,16.00){\makebox(0,0)[cc]{$_{4_d}$}}
\put(137.00,16.00){\makebox(0,0)[cc]{$_{4_a}$}}
\put(167.00,16.00){\makebox(0,0)[cc]{$_{4_e}$}}
\put(197.00,16.00){\makebox(0,0)[cc]{$_{4_d}$}}
\put(112.00,27.00){\makebox(0,0)[cc]{$_2$}}
\put(112.00,32.00){\makebox(0,0)[cc]{$_5$}}
\put(142.00,27.00){\makebox(0,0)[cc]{$_2$}}
\put(142.00,32.00){\makebox(0,0)[cc]{$_5$}}
\put(135.00,30.00){\makebox(0,0)[cc]{$_3$}}
\put(133.00,26.00){\makebox(0,0)[cc]{$_5$}}
\put(133.00,32.00){\makebox(0,0)[cc]{$_2$}}
\put(172.00,27.00){\makebox(0,0)[cc]{$_2$}}
\put(172.00,32.00){\makebox(0,0)[cc]{$_5$}}
\put(165.00,30.00){\makebox(0,0)[cc]{$_3$}}
\put(163.00,26.00){\makebox(0,0)[cc]{$_5$}}
\put(163.00,32.00){\makebox(0,0)[cc]{$_2$}}
\put(195.00,30.00){\makebox(0,0)[cc]{$_3$}}
\put(193.00,26.00){\makebox(0,0)[cc]{$_5$}}
\put(193.00,32.00){\makebox(0,0)[cc]{$_2$}}
\put(120.00,26.00){\makebox(0,0)[cc]{$_{0_d}$}}
\put(126.00,26.00){\makebox(0,0)[cc]{$_{0_c}$}}
\put(150.00,26.00){\makebox(0,0)[cc]{$_{0_a}$}}
\put(156.00,26.00){\makebox(0,0)[cc]{$_{0_f}$}}
\put(180.00,26.00){\makebox(0,0)[cc]{$_{0_e}$}}
\put(186.00,26.00){\makebox(0,0)[cc]{$_{0_b}$}}
\put(122.00,37.00){\makebox(0,0)[cc]{$_4$}}
\put(122.00,21.00){\makebox(0,0)[cc]{$_4$}}
\put(152.00,37.00){\makebox(0,0)[cc]{$_4$}}
\put(152.00,21.00){\makebox(0,0)[cc]{$_4$}}
\put(182.00,37.00){\makebox(0,0)[cc]{$_4$}}
\put(182.00,21.00){\makebox(0,0)[cc]{$_4$}}
\put(109.00,28.00){\makebox(0,0)[cc]{$_3$}}
\put(139.00,28.00){\makebox(0,0)[cc]{$_3$}}
\put(169.00,28.00){\makebox(0,0)[cc]{$_3$}}
\put(44.00,45.00){\vector(-1,1){8.00}}
\put(34.00,53.00){\vector(-1,-1){8.00}}
\put(26.00,43.00){\vector(1,-1){8.00}}
\put(36.00,35.00){\vector(1,1){8.00}}
\put(66.00,35.00){\vector(1,1){8.00}}
\put(74.00,45.00){\vector(-1,1){8.00}}
\put(64.00,53.00){\vector(-1,-1){8.00}}
\put(56.00,43.00){\vector(1,-1){8.00}}
\put(14.00,45.00){\vector(-1,1){8.00}}
\put(94.00,53.00){\vector(-1,-1){8.00}}
\put(86.00,43.00){\vector(1,-1){8.00}}
\put(6.00,35.00){\vector(1,1){8.00}}
\put(6.00,5.00){\vector(1,1){8.00}}
\put(14.00,15.00){\vector(-1,1){8.00}}
\put(36.00,5.00){\vector(1,1){8.00}}
\put(44.00,15.00){\vector(-1,1){8.00}}
\put(34.00,23.00){\vector(-1,-1){8.00}}
\put(26.00,13.00){\vector(1,-1){8.00}}
\put(66.00,5.00){\vector(1,1){8.00}}
\put(74.00,15.00){\vector(-1,1){8.00}}
\put(64.00,23.00){\vector(-1,-1){8.00}}
\put(56.00,13.00){\vector(1,-1){8.00}}
\put(94.00,23.00){\vector(-1,-1){8.00}}
\put(86.00,13.00){\vector(1,-1){8.00}}
\put(94.00,24.00){\vector(-1,0){28.00}}
\put(146.00,30.00){\vector(-1,1){8.00}}
\put(136.00,38.00){\vector(-1,-1){8.00}}
\put(128.00,28.00){\vector(1,-1){8.00}}
\put(138.00,20.00){\vector(1,1){8.00}}
\put(116.00,30.00){\vector(-1,1){8.00}}
\put(108.00,20.00){\vector(1,1){8.00}}
\put(196.00,39.00){\vector(-1,0){28.00}}
\emline{168.00}{19.00}{33}{196.00}{19.00}{34}
\put(35.00,54.00){\circle{2.00}}
\put(64.00,24.00){\vector(-1,0){28.00}}
\put(34.00,24.00){\vector(-1,0){28.00}}
\put(94.00,54.00){\vector(-1,0){28.00}}
\put(64.00,54.00){\vector(-1,0){28.00}}
\put(34.00,54.00){\vector(-1,0){28.00}}
\put(6.00,4.00){\vector(1,0){28.00}}
\put(36.00,4.00){\vector(1,0){28.00}}
\put(66.00,4.00){\vector(1,0){28.00}}
\put(166.00,39.00){\vector(-1,0){28.00}}
\put(136.00,39.00){\vector(-1,0){28.00}}
\put(108.00,19.00){\vector(1,0){28.00}}
\put(138.00,19.00){\vector(1,0){28.00}}
\put(6.00,34.00){\vector(1,0){28.00}}
\put(36.00,34.00){\vector(1,0){28.00}}
\put(66.00,34.00){\vector(1,0){28.00}}
\put(196.00,39.00){\vector(-1,0){28.00}}
\put(168.00,19.00){\vector(1,0){28.00}}
\put(168.00,20.00){\vector(1,1){8.00}}
\put(176.00,30.00){\vector(-1,1){8.00}}
\put(166.00,38.00){\vector(-1,-1){8.00}}
\put(158.00,28.00){\vector(1,-1){8.00}}
\put(196.00,38.00){\vector(-1,-1){8.00}}
\put(188.00,28.00){\vector(1,-1){8.00}}
\put(197.00,29.00){\vector(0,-1){9.00}}
\put(197.00,29.00){\vector(0,1){9.00}}
\put(167.00,29.00){\vector(0,-1){9.00}}
\put(167.00,29.00){\vector(0,1){9.00}}
\put(137.00,29.00){\vector(0,-1){9.00}}
\put(137.00,29.00){\vector(0,1){9.00}}
\put(107.00,29.00){\vector(0,-1){9.00}}
\put(107.00,29.00){\vector(0,1){9.00}}
\put(95.00,44.00){\vector(0,-1){9.00}}
\put(95.00,44.00){\vector(0,1){9.00}}
\put(65.00,44.00){\vector(0,-1){9.00}}
\put(65.00,44.00){\vector(0,1){9.00}}
\put(95.00,14.00){\vector(0,-1){9.00}}
\put(95.00,14.00){\vector(0,1){9.00}}
\put(65.00,14.00){\vector(0,-1){9.00}}
\put(65.00,14.00){\vector(0,1){9.00}}
\put(35.00,44.00){\vector(0,-1){9.00}}
\put(35.00,44.00){\vector(0,1){9.00}}
\put(5.00,44.00){\vector(0,-1){9.00}}
\put(5.00,44.00){\vector(0,1){9.00}}
\put(35.00,14.00){\vector(0,-1){9.00}}
\put(35.00,14.00){\vector(0,1){9.00}}
\put(5.00,14.00){\vector(0,-1){9.00}}
\put(5.00,14.00){\vector(0,1){9.00}}
\end{picture}
\caption{Split representation of the quotient graph $D/\Z_7$}
\end{figure}
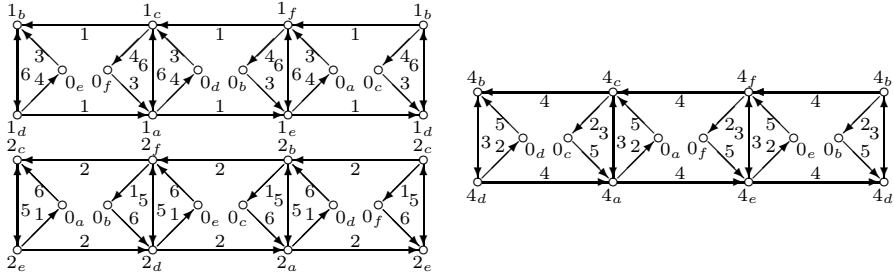

With these symbols adopted, the quotient graph $D/\Z_7$ can be considered as a voltage graph with group $\Z_7$ and derived graph $D$, (\cite{GT}), admitting a split representation into the three connected digraphs of Figure 2, whose vertices are indicated by the symbols $j_i$ of their representatives $yup_0$, and in which:
{\bf(1)}
the 18 oriented 4-cycles that are shown are interpreted all with counterclockwise orientation;
{\bf(2)} For each $i\in\{a,\ldots,f\}$, the three vertices indicated by $0_i$ represent just
one vertex of $D/Z_7$, so they must be identified;
{\bf(3)} the leftmost arc in each one of the three connected digraphs must be identified with
the corresponding rightmost arc by parallel translation;
{\bf(4)} if an arc $\vec{e}$ of $D/\Z_7$ has voltage $\nu\in\Z_7$, initial vertex $j_i$ and
terminal vertex $j'_{i'}$, then a representative
$yup_\mu$ of $j_i$ initiates an arc in $D$ that covers $\vec{e}$ and has terminal vertex $y'u'p'_{(\nu+\mu)}$, where $yup_0$ and $y'u'p'_0$ are represented respectively by $j_i$ and $j'_{i'}$.

All the oriented 4-cycles of $D$ are obtained by uniform translations mod 7 from these 18 oriented 4-cycles. Thus, there are just $126=7\times 18$ oriented 4-cycles of $D$.
Our construction of $D$ shows that the following statement holds.

\begin{thm} The oriented graph $D$ is a strongly connected $\vec{C}_4$-UH digraph on $168$ vertices, $126$ pairwise disjoint oriented $4$-cycles, with regular indegree and outdegree both equal to $3$ and no circuits of lengths $2$ and $3$.\qfd\end{thm}

\end{document}